# A Efficient Algorithm For Computing k-Factor

Yingtai  Xie
Chengdu University

**Abstract**

This article will prove a theorem for the existence of *k*-factor for $k>1$, and present an efficient algorithm for computing *k*-factor for all values of $k \geq 2$ based on this theorem.

## 0.Introduction

In this article, both of graph G and *k*-factor of G are simple, i.e there are no multiple edges (connecting the same pair of vertices) and no loops. A *k*-regular spanning subgraph of *G* is called a *k*-factor of G. A perfect matching of G is a 1-factor of G. A spanning subgraph of G, which is a set of disjoint cycles, is a 2-factor of G. The *k*-factor, especially 2-factor, has many applications in Graph Theory, Computer Graphics, and Computational Geometry. Umans and Lenhart [2] constructed an algorithm for solid grid graph to find a Hamltonian cycle by merging cycles of a 2-factor. Yingtai Xie [11][12] use 2-factor in algorithm for fingding a Hamilton cycle/path in 2D-mesh, which is the topologies for Networks on Chip and on huge number of mesh-connected computer systems. In Computer Graphics, 2-factors be used in algorithms for fast rendering of 3D scenes [4]. Similarly, 2-factors be also used for stripification of large sets of quadrilaterals and tetrahedra [3].

1-factor is a *perfect matching* of a graph $G$. Matching problem, which occupies a central place in graph algorithm theory, has been thoroughly studied. Based on the *M*-augmenting path ideal introduced in well known Berge's theorem [1], many effective algorithms for computing 1-factore are established. However, Berge's theorem does not apply to *k*-factor, when $k \geq 2$, therefor the algorithm for computing 1-factor based on M-augmenting path can not be applied to compute *k*-factor of general graph directly.

In the section 2 of this article, we introduce the ***M-augmenting trail*** ideal and prove a theorem for the existence of *k*-factor for $k>1$, and an algorithm for computing k-factor for k>1 based on this theorem is presented In the section 3.

## 1. Previous work

Petersen (1891) proved that every regular graph of positive even degree has a 2-factor. Petersen's result establishes the existence of 2*k*-factors in 2*m*-regular graphs, $m \geq k$, and his proof leads a algorithm for computing 2*k*-factor of a 2*m*-regular graph, which only need to compute 1-factor in a bipartite *m*-regular graph. This is the first time that the algorithm for computing 1-factor was applied to compute *k*-factor for a class of graph. Umans and Lenhart [2] present a algorithm for computing 2-fctor in solid grid graph. Yingtai Xie [11][12] provide a algorithm based on Augmenting path to computing 2-factor in general grid graph.

Until now, the mainstream algorithm for computing *k*-factors in general graphs is to derive a (transformed) graph from the original graph such that the original graph has a k- factor, if and only if the derived graph has a 1-factor. The method is to design a template for each point with degree *d* of the input graph G, and replace each point with its template to derive a graph G'. Thus, a *k*-factor of the original graph can be found by computing the perfect matching of the derive graph. Diaz-Gutierrez and Gopi [4] for vertices of degree 4 and Umans [5] for grid graph designed template substitution. H. Meijer et al [6] generalize to vertices of any degree and to any possible *k* value designed template substitution for computing *k*-factor .

A template has more vertices , for example, the template designed by Diaz-Gutierrez and Gopi [4] and Uman[5] for the vertex of degree 4 has six vertices and number of points of a template increases with the degree of the vertices of original graph ,H. Meijer et al [6] denote that let *G* be the input graph, $|V(G)|=n$ , and *G'* be the transformed graph obtained from *G* through the template substitution then $|V(G')|=o(n^2)$ , $|E(G')|=o(n^2)$ the time complexity for computing 1-factor is $o(|V(G')|^2 \sqrt{|E(G')|})=o(n^3)$ , by Micali and Vazirani [7,8] , in transformed graph *G'* , therefor the time complexity for computing *k*-factor is $o(n^3 k)$ in input graph G with $|V(G)|=n$ ,the timing of the template transformation is not taken into account.

## 2.The theorem of existence for k-factor.

We call *M* **k-limited spanning subgraph** of *G* if $V(M)=V(G)$ and $d_M(v) \le k$ for each $v \in V(G)$ . Vertex *v* is called **filled** if $d_M(v)=k$ ,otherwise **unfilled**. $E(M)$ and $E(G-M)$ are called red and blue respectively .Let

$$\sigma(M) = \sum_{v \in V(G)} (k - d_M(v)) \qquad (2.1)$$

$\sigma(M) \ge 0$ is even and *M* is *k*-factor of *G* if $\sigma(M)=0$ .

A **trail** (of length *m*) in a graph *G* is a non-empty sequence $P = v_0 e_1 v_1 e_2 v_2,...v_{m-1} e_m v_m$ of vertices and distinct edges in *G* such that $e_i=(v_{i-1},v_i)$ for all $i \le m$ . If $v_0 = v_m$, the trail is closed.

A trail *P* is called **M-alternating trail** if its edges are alternately of blue and red .

**Definition 2.1**  A *M- augmenting trail*

$$P = v_0 e_1 v_1 e_2 v_2,...v_{m-1} e_m v_{2m} \qquad (2.2)$$

*is a M* -alternating trail such that its end points $v_0$ and $v_{2m}$ both are unfilled in *M* and both of the first edge $e_1$ and the terminal edge $e_{2m}$ are blue. And all inner vertex $v_1 v_2,...,v_{2m-1}$ of *P* are filled in *M* .and $d_M(v_0) < k-1$ is required for $v_0 = v_{2m}$ .

If *P* is a M-augmenting trail then $|P|$ must be odd .

**Lemma 2.1**  let *M* be a *k-limited* spanning subgraph, $\sigma(M) \ne 0$ , and *P* is a M-augmenting trail with end points $v_0$ and $v_k$ then

$$M' = M \oplus P \qquad (2.3)$$

Is also a *k-limited* spanning subgraph of *G* and $\sigma(M') = \sigma(M) - 2$ .

*proof :*  $M'$ is constructed by deleting all the red edges of *P* from *M* and adding to *M* all blue edges of *P* . It is clear that $d_{M'}(v_i) = d_M(v_i)$ for all of the filled vertices $v_i, i=1,2,,..2m-1$ ,because the number of blue edges is the same as the number of red edges at the filled vertices and just are each added a blue edge at the unfilled ends *v₀* and *v₂ₘ,* hence $d_{M'}(v_0) = d_M(v_0)+1, d_{M'}(v_{2m}) = d_M(v_{2m})+1$ , therefore $\sigma(M_1) = \sigma(M) - 2$ . (refer to Fig.3.1,(a) A path made up of red edges is a 2-limited spanning subgraph M with two unfilled vertices 13 and 18. $P = (13,5,6,12,11,9,10,6,7,18,16,19)$ is a M-augmenting trail. (b) $M' = M \oplus P$ is a 2-factor.)

From this, we can infer directly.

**Theorem 2.1**  $G$ have a $k$-factor if and only if there is a $M$-augmenting trail for any $k$-limited spanning subgraph $M$ with $\sigma(M) > 0$.

*Proof* *'if'* can be infer from lemma 2.3 directly.

If $G$ have a $k$-factor $F$. Consider the graph $N = F \oplus M$, Since $d_F(v) = k$ and $d_M(v) \leq k$ then $d_N(v) \leq 2k$ for each vertex $v \in V(G)$.

Let edge $e$ be red if $e$ is of $F$, and blue if $e$ is of $M$.

(i) If $v$ is filled vertex of $M$, then $d_F(v) = d_M(v) = k$

  (i$_a$) $d_N(v) = 0$, $v$ is a isolated in N.

  (i$_b$) $0 < d_N(v) \leq 2k$ is even and there are even edges at $v$, half of them are red and another half are blue.

(ii) If $v$ is unfilled of $M$, i.e $d_M(v) < k$, then

There are $s + t$ edges at $v$, $s$ blue edges and $t$ red edges, $s < t$, $s + t < 2k$.

since $\sigma(M) > 0$, hence the number of blue edges is less than the number of red edges, there's at least one connected component N' of N , which either include at least two unfilled vertices of $M$ or only one unfilled vertex of dgree< $k$-1 in $M$ .Let $v_0$ of N' be a unfilled vertex, $v_0$ must be incident with a red edge $e$, if $e$ is incident another unfilled vertex $u$, then $v_0 e u$ is a M-augmenting trail,otherwise follow a trail of edges red and blue alternating from $u$ until pass a blue edge to a unfilled vertex $v_m$.

$$P = v_0 e_1 v_1 e_2 v_2 \ldots e_m v_m \qquad (2.4)$$

It is not possible to get stuck at any filled vertex, because when the trail enters a filled vertex there must be an unused alternating edge leaving it by (i$_b$) ,as the same reason,it is not possible to get stuck at any blue edge by (ii). $d_M(v_0) < k - 1$ is required for $v_0 = v_m$. The trail (2.3) formed in this way is just a M-augmenting trail. □

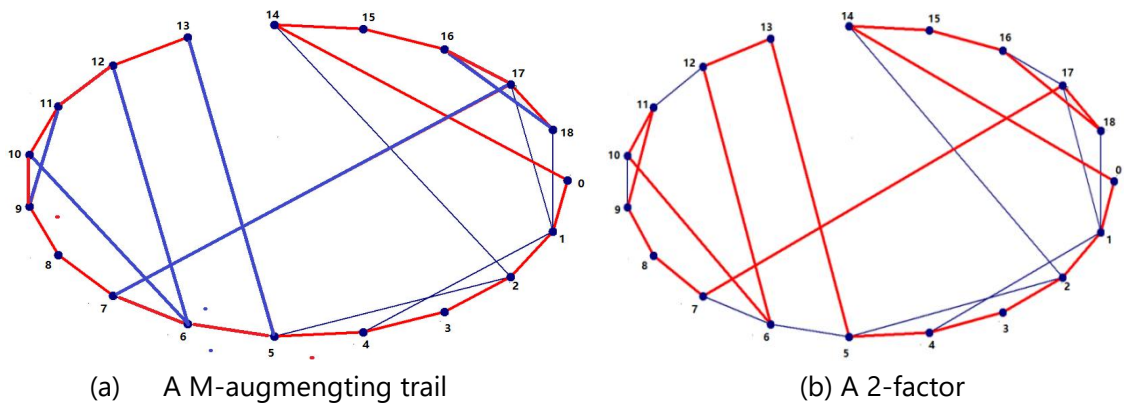

(a)  A M-augmengting trail      (b) A 2-factor

Fig. 2.1 A illustration for Lemma 2.1 and Theorem 2.1

## 3.Algorithm for k-factor

Algorithm A:algorithm for k-factor

$$M = \phi$$

do

    Find a *M*-augmenting trail *P*.

    If there is no *P* then there is no k-factor;exit do

      else

$$M := M \oplus P$$

    If $\sigma(M) = 0$ then *M* is a k-factor; exit do

loop

Algorithm B.  Algorithm for M-augmenting trail.

The M-augmenting trail *P* with minimal cardinality is called shortest. It is clear that there is no even alternating cycles in a shortest M-augmenting trail ,therefor

**Claim 3.1** A shortest M-augmenting trail (Eq.(2.2)) in a bipartite graph is just a M-augmenting path.

Therefore, it can be deduced directly from Theorem 2.1 that,

**Theorem 3.1** The algorithm based on Augmenting path can be used to computing k-factor of the bipartite graph.

We call an odd alternating cycle a ***blossom.***

Different from the augmenting path in the matching algorithm(1-factor) in general graph, it is possible that there are odd cycles(blossom) in a M-augmenting trail when $k \geq 2$ .therefor Edmond's method of dealing with blossoms by "shrinking" them[9][10]，which is the basic in the algorithm for computing 1-factor of the general graph, are not adequate for the task of finding M-augmenting trail，since the blossoms are completely lost.

To construct the algorithm for M-augmenting trail,we associate with each edge $e = (u,v)$ of $G$ two distinct elements $d = \overrightarrow{uv}$ and $d^- = \overrightarrow{vu}$ called the opposite darts on *e* and associate each dart with two end points of $d$ ,one is called its head $h(d)$ , other is called its tail $t(d)$ . We construct a equivalent directed graph $\hat{G} = (\hat{E},\hat{V})$ of $G$ such that $\hat{V} = V(G)$ and $\hat{E} = \{d, d^- : d \text{ and } d^- \text{ are associed with } e \in E(G)\}$ . And  $d \in \hat{M}$  or  $d \in \hat{E} - \hat{M}$ according $e \in M$ or $e \in E - M$ .

A ***directed trail*** in $\hat{G}$ is a sequence

$$\vec{P} = d_0 d_1, \ldots d_k \qquad (3.1)$$

of $k+1$ distinct darts of $\hat{G}$ .It is required that that the head of $d_j$ shall be the tail of $d_{j+1.}$

A ***directed trail*** in $\hat{G}$ is *edge-simple* if no two of its terms are darts on the same edge of $G$, i.e. there is no $d_i = d_j^-$ in $\vec{P}$.

**Definition 3. 1**  A directed trail $\vec{P}$ is called ***M-alternating directed trail*** if its darts are alternately of $\hat{E} - \hat{M}$ and $\hat{M}$ .

**Definition 3.2**  A  *M*-alternating directed edge-simple trail $\vec{P} = d_0 d_1, \ldots d_{2m}$ is called ***M-***

***augmenting directed trail*** if its end points $v_0 = t(d_0)$ and $v_{2m} = h(d_{2m})$ both are unfilled in $M$ and both of the first dart $d_0$ and the terminal $d_{2m}$ are of $\vec{E} - \vec{M}$. And all inner vertex $v_i = h(d_i)$ $0 \leq i \leq 2m-1$ of $\vec{P}$ are filled in $M$.

**Claim 3.2** It is clear that if $\vec{P} = d_0 d_1, \ldots d_{2m}$ is a M-augmenting directed trail then $P = v_0 e_1 v_1 e_2 . e_{2m} v_{2m}$ is a M-Augmenting trail, where $v_0 = t(d_0)$, $v_i = h(d_i)$ $0 \leq i \leq 2m$, And $e_i$ is edge which the dart $d_i$ is on.

To construct an algorithm for finding a M-augmenting directed trail in $\hat{G}$, we construct a directed graph $\vec{G}$.

let $D$ be a set of darts, $t(D) = \{t(d) \mid d \in D\}$ and $h(D) = \{h(d) \mid d \in D\}$. To construct a sequence of $D_i$ $i=0,1,2,k,$ definition:

a) $t(D_0) = \{u : d_M(u) < k\}$
b) $t(D_{i+1}) \subseteq h(D_i)$ for i=0,1,2,...
c) $D_i \subset \vec{E} - \vec{M}$ for $i = 0, 2, 4, \ldots$, $D_i \subset \vec{M}$ for $i = 1, 3, \ldots$,

until for the first time a set $D_{2m}$ is constructed such that $h(D_{2m})$ include an unfilled vertex $v$ of $M$( $d_M(u) < k - 1$ is required for $v = u$ ).

Let

$$\vec{G} = D_0 \cup D_1 \cup \ldots \cup D_{2m} \qquad (3.2)$$

(refer to Fig.3.1 (a) the directed graph starting at a unfilled vertex $u=13$ of fig3.1 (a))

In order to find directed $M$ - augmenting directed trail starting at $u$ and ending to $v$, it is necessary to delete such vertices which cannot be on a directed M- augmenting directed trail, and the deletion of such vertices may make certain other vertices breakpoints, in the sense that they are not heads or tails of any darts. These breakpoints must also be deleted, hence we need an algorithm which, deletes some vertices of $\vec{G}$, together with all other vertices maked breakpoints by this deletion.

This is achieved as follows:

We call vertex $x$ dead-Head if $x \in h(D_i)$ and there is no $d \in D_{i+1}$ with $x = t(d)$, and $x$ dead-Tail if $x \in t(D_i)$ and there is no $d \in D_{i-1}$ with $x = h(d)$

***Pruning Algorithm:***

1. delete dead-Head vertices from $h(D_i)$ in order $i = 2m, 2m-1, \ldots, 3, 2, 1, 0$
2. delete dead-T vertices from $t(D_i)$ in order $i = 0, 1, 2, \ldots 2m$

We just have to find directed $M$-augmenting trail in the subgraph $\vec{G}_v$ of $\vec{G}$ which contains and only contains all of the alternating-trails ending to $v$ from $u$.

***Algorithm for constructing for*** $\vec{G}_v$: Delete all vertices except $v$ from $h(D_{2m})$ and then do ***Pruning Algorithm*** for result graph. (refer to Fig.3.1 (b) the directed graph $\vec{G}_{18}$ starting at $u=13$ and ending to $v=18$ of Fig.3.1 (a) )

**Claim 3.1**. $\vec{G}_v$ is set of the directed alternating trails starting at $u$ and ending to $v$.

***Proof*** :If $\vec{P} = d_0 d_1, \ldots d_{2m}$ of $\vec{G}_v$ have repeating darts then there will be a walk ending to

$v$, which is shorter than $\vec{P}$, this contradicts the fact that $v \in h(D_{2m})$ only (3.2), therefor $\vec{P}$ can only be a direct trail. □

An alternating-trail $\vec{P} = d_0 d_1, \ldots d_{2m}$ pass a blossom and is not egde-simple then there must be two darts $d_i = d_j^-$, $i < j, h(d_i) = t(d_j)$. We call $d_i$ **in-dart** and $d_j$ **out-dart** of the blossom.

We can find out a directed M-augmenting path in $\vec{G_v}$ by **Blossom** operation.

**Blossom operation:** if the in-dart of a blossom is not a cut-dart of $\vec{G_v}$ then delete it else delete out-dart from $\vec{G_v}$.

### The algorithm for construct an M-augmenting trail.

a) Construct the directed graph $\vec{G}$ .(3.2)

b) Construct the directed graph $\vec{G_v}$ from $\vec{G}$.

c)
   Do
       Get a M-alternating-trail $\vec{P} = d_0 d_1, \ldots d_{2k}$ from $\vec{G_v}$.
         If $\vec{P}$ is edge-simple then to get a M-augmenting trail by claim 3.2.
       Else
         Do blossom operation on $\vec{G_v}$
   loop

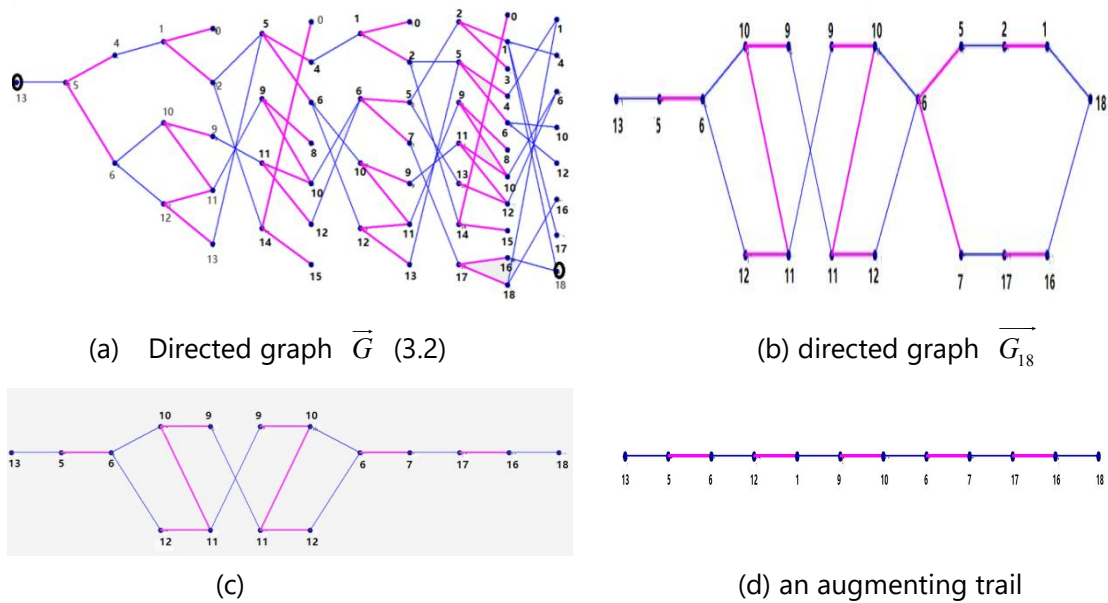

(a) Directed graph $\vec{G}$ (3.2)        (b) directed graph $\vec{G_{18}}$

(c)        (d) an augmenting trail

Fig.3.1 An example of an application of the algorithm for construct an M-augmenting trail.

**Example** 3.1. The graph $G$ shown by Fig.1.1 (a) have a 2-limited subgraph $M$ with two unfilled end-vertices 13 and 18 (a path from 13 to 18 constructed by red edges ).The Directed graph

$$\vec{G} = D_0 \cup D_1 \cup \ldots \cup D_{10}$$

Is shown by Fig 3.1 (a), $t(D_0) = \{13\}$ and $\{18\} \subset h(D_{10})$. The directed graph $\vec{G_{18}}$ is constructed by Deleting all vertices except 18 from $h(D_{10})$ and then do **Pruning Algorithm** for result graph, which as shown by Fig 3.1 (b). $\vec{P} = (13, 5, 6, 10, 9, 11, 10, 6, 5, 2, 1, 18)$ is a directed trail from 13 to 18 of $\vec{G}_{18}$, which pass a blossom (6,10,9,11,12,6) and $d_1 = \overrightarrow{(5,6)}$ is in-dart, $d_7 = \overrightarrow{(6,5)}$ is out-dart of it. $d_1 = \overrightarrow{(5,6)}$ is cut-dart of $\vec{G}_{18}$. so doing blossom operation on $\vec{G}_{18}$ to get rid of the out-dart $d_7 = \overrightarrow{(6,5)}$. The result graph is shown by Fig.3.1 (c), to get a directed trail $\vec{P} = \overrightarrow{(13,5,6,10,9,11,10,6,7,17,16,18)}$ in it, which pass a blossom $\overrightarrow{(10, 9, 11, 10)}$ with in-dart $d_3 = \overrightarrow{(6,10)}$ and out-dart $d_6 = \overrightarrow{(10,6)}$, the in-dart is not a cut-dart, to get rid of it, the result graph is a augmenting trail, as shown by Fig.3.1(d).

The algorithm for constructing an M-augmenting trail is to construct $\vec{G}$, and $\vec{G}_v$ first, and then followed by an iterative search on $\vec{G}_v$: To find a directed alternating trail, if it is not a augmenting trail then executes the blossom operation on $\vec{G}_v$. The iteration stops when an augmenting trail have been found. It is easy to see that each edge of G is examined once when $\vec{G}$ is first constructed, and once more if an edge or an end point of that edge is deleted when doing Pruning Algorithm.

It is easy to see that the execution time of a search is $o(m)$, where $m$ is the number of edges in G, if $n$ is the number of vertices, then $k(n/2)$ searches at most, hence the execution time of the entire algorithm is $o(kmn)$.

## 4. Conclusion

The idea of Augmenting path was first introduced by Berge (1957). Today, it has become a very mature and standard technique in matching theory. Berger's theorem gives the characteristic of maximum matching, based on which many efficient algorithms for computing 1-factor are constructed.

This article introduces the idea of augmenting trail, which will play the same role in the theory and algorithm of K-factor (k>1) as augmenting path in the study of 1-factor. Based on the augmenting trail, an efficient algorithm for computing k-factor are constructed, in addition this algorithm deals directly and simply with the blossom.

which have been plagued the formation of the 1-factor algorithm for many years